\documentclass[11pt,twoside]{article}

\usepackage[ansinew]{inputenc}
\usepackage{amsfonts}
\usepackage{amssymb,amsmath}
\usepackage{latexsym}

\newcommand {\debeq}	{\begin{eqnarray*}}
\newcommand {\fineq}	{\end{eqnarray*}}

\newcommand	{\tendinfty}	
{\rightarrow\infty}

\newcommand	{\ZZ}{\mathbb{Z}}

\newcommand	{\NN}{\mathbb{N}}

\newcommand	{\RR}{\mathbb{R}}
\newcommand	{\PP}{\mathbb{P}}
\newcommand	{\EE}{\mathbb{E}}

\newtheorem	{thm}		{Theorem}[section]

\newtheorem     {rem}           {Remark}
\newtheorem	{prop}	[thm]{Proposition}

\begin{document}
\title{
Past, growth and persistence of source-sink metapopulations
}
\author{\textsc{By Vincent Bansaye and  Amaury Lambert}
}
\date{}								
\maketitle
\noindent

\medskip

\noindent
\textsc{Vincent Bansaye\\
CMAP\\
\'Ecole Polytechnique\\
Route de Saclay\\
F-91128 Palaiseau Cedex, France}\\
\textsc{E-mail: }vincent.bansaye@polytechnique.edu\\
\textsc{URL: }http://www.cmapx.polytechnique.fr/$\sim$bansaye/
\\
\\
\textsc{Amaury Lambert\\
Laboratoire de Probabilités et Modèles Aléatoires\\
UMR 7599 CNRS and UPMC Univ Paris 06\\
Case courrier 188\\
4, Place Jussieu\\
F-75252 Paris Cedex 05, France}\\
\textsc{E-mail: }amaury.lambert@upmc.fr\\
\textsc{URL: }http://www.proba.jussieu.fr/pageperso/amaury/index.htm
\\
\\

\begin{abstract}
\noindent

Source-sink systems are metapopulations of patches that can be of variable habitat quality.  They can be seen as graphs, where  vertices represent the patches, and the weighted oriented edges give the probability of dispersal from one patch to another. We consider either finite or source-transitive graphs, i.e., graphs that are identical when viewed from a(ny) source. We assume stochastic, individual-based, density-independent reproduction and dispersal.

By studying the path of a single random disperser, we are able to display simple criteria for persistence, either necessary and sufficient, or just sufficient. In case of persistence, we characterize the growth rate of the population as well as the asymptotic occupancy frequencies of the line of ascent of a random survivor.  
Our method allows to decouple the roles of reproduction and dispersal.
%
%
%
%
Finally, we extend our results to the case of periodic or random environments, where some habitats can have variable growth rates, autocorrelated in space and possibly in time. 

In the whole manuscript, special attention is given to the example of regular graphs where each pair of adjacent sources is separated by the same number of identical sinks. In the case of a periodic and random environment, we also display examples where all patches are sinks when forbidding dispersal but the metapopulation survives with positive probability in the presence of dispersal, as previously known for a two-patch mean-field model with parent-independent dispersal \cite{JansYosh}.

\end{abstract}  	

\tableofcontents

\section{Introduction}

\paragraph{Model.}
We consider a stochastic, individual-based model of spatially structured population dynamics. 
The spatial structure is a metapopulation of patches that can be of different habitat qualities. We label by $i=1,\ldots,K$ the patches so that the model can  be described by a labeled finite graph with weighted oriented edges. Vertices represent the patches, an oriented edge from vertex $i$ to vertex $j$ bears a weight $d_{ij}$ equal to the probability of dispersal from patch  $i$ to patch $j$.

We assume a simple asexual life cycle with discrete non-overlapping generations and no density-dependence. At each generation, as a result of survival and reproduction, all individuals, independently from one another, leave to the next generation a random number of individuals, called offspring, whose mean number depends on the habitat quality. Immediately after local growth, each individual from the new generation migrates independently, from patch $i$  to patch $j$ with probability $d_{ij}$. No mortality is assumed during dispersal, since it can be encompassed in the growth phase. Reproduction, survival and dispersal probabilities are assumed not to depend on local densities. 

The mean \emph{per capita} number of offspring in patch $i$ will be denoted by $m_i$, and can be seen as a proxy for habitat quality. If $m_i>1$, the population may persist without dispersal and we say that patch $i$ is a \emph{source}. Conversely, if $m_i\leq 1$, the population would die out without dispersal and  we say that patch $i$ is a \emph{sink}.  It will be convenient to assume that $m_1\geq m_2\geq \cdots\geq m_K$ and more relevant to consider the case when $m_1>1\geq m_K$.  Even in the presence of sinks, the metapopulation might persist thanks to local growth on sources replenishing sinks by dispersal. 
  In the second part of the paper, we will also assume that environment is variable and denote by $m_i(w)$ the mean offspring in a patch of type $i$ when the environment is in state $w$. 
    Finally, we extend our results to source-sink metapopulations on transitive graphs. We also want to make the observation that our approach applies to more general multitype branching processes, in the sense that dispersal behaviours of siblings might not be independent. 

\begin{center}
\begin{figure}[ht]
\unitlength 1mm 
\linethickness{0.4pt}
\input{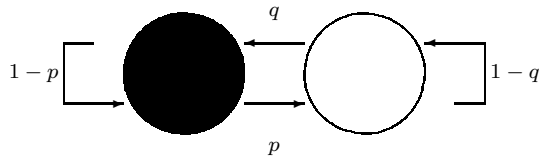}
\label{fig:two-patches}
\caption{Two patches of different qualities. The filled circle is a source and the empty circle is a sink. The arrow labels are the dispersal probabilities.}
\end{figure}
 \end{center}

\begin{center}
\begin{figure}[ht]
\unitlength 1.3mm 
\linethickness{0.4pt}
\input{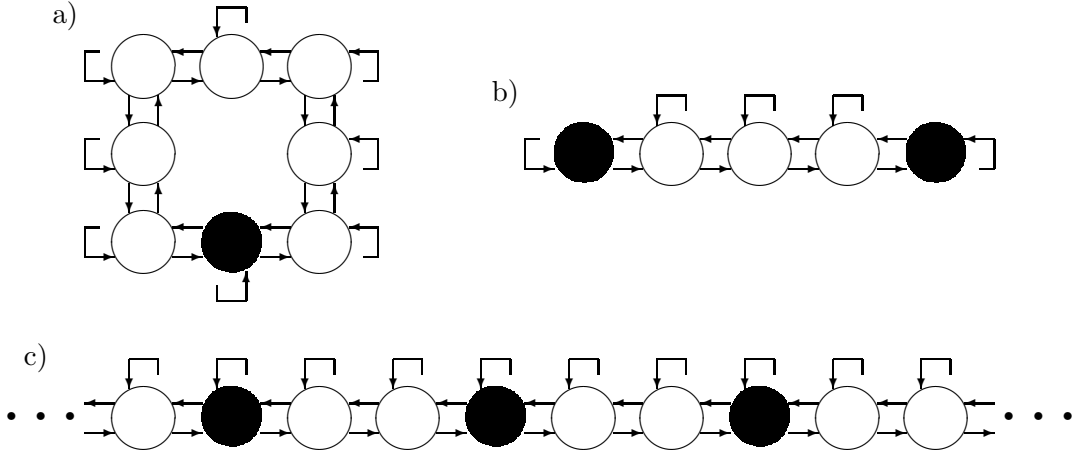}
\caption{Three examples where sources are only connected to sinks and two adjacent sources are connected by $n$ sinks; a) cyclic graph with one source, $n=7$; b) two sources, $n=3$; c) infinite periodic array, $n=2$, arrow labels (not represented) are also assumed periodic.}
\label{fig:arrays}
\end{figure}
 \end{center}

\paragraph{Two natural examples.} Let us describe two examples with two possible habitat qualities, one source type and one sink type, that will be treated throughout the paper:
$$M:=m_1>1, \qquad m:=m_2=\cdots=m_K\leq 1.$$ 
 First, we are interested in the simple case with two patches, patch $1$ with mean offspring  $M$, and patch $2$ with mean offspring  $m$. In this case, we will always use the simplified notation $p=d_{12}$ and $q=d_{21}$ (see Figure \ref{fig:two-patches}). \\
Second, we will consider the case  when each source is only connected to sinks and two adjacent sources are separated by a line of $n$ identical sinks. An example of such graph is the cyclic finite graph with one source and $n$ sinks, or two sources connected by $n$ sinks, or an infinite array with period $n$ (see Figure \ref{fig:arrays})...
 
\paragraph{Outline of the paper.}
A random disperser is   a  walker on the graph which follows the  dispersal stochastic scheme. We display here  a very simple persistence criterion involving the mean growth rates in each type of habitat and the times spent in each habitat, between two visits of a source, by a single random disperser.  We specialize this result for the examples given above. 
  We also give the asymptotic fraction of time spent in each patch by an individual taken at random in the surviving population. We find that this distribution never equals the stationary distribution of the single disperser, except in the case when all habitat qualities are identical.

Addressing those questions is much more difficult when the environment is variable, i.e., when each habitat type can have a different mean growth rate depending on the state of the environment. Applying our techniques yields partial results, viz. when the environment is periodic or given by an ergodic sequence of random variables. In these two cases, we prove that the population can survive in sinks only, in the sense that without dispersing, it would die out in each patch, a result which was previously known \cite{JansYosh} in the case of large populations on an island model (parent-independent dispersal). 

\begin{center}
\begin{figure}[ht]
\unitlength 1mm 
\linethickness{0.4pt}
\input{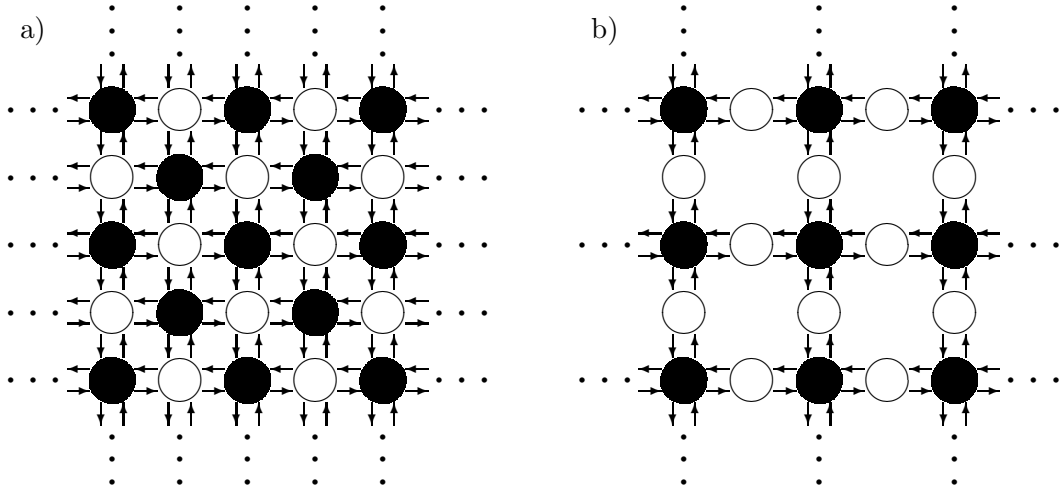}
\caption{Two examples of source-transitive graphs with a finite motif  ; a) the chessboard; b) a square grid where four-degree vertices are sources separated by $n$ sinks (here $n=1$).}
\label{fig:transitive}
\end{figure}
\end{center}

\begin{center}
\begin{figure}[ht]
\unitlength 1mm 
\linethickness{0.4pt}
\input{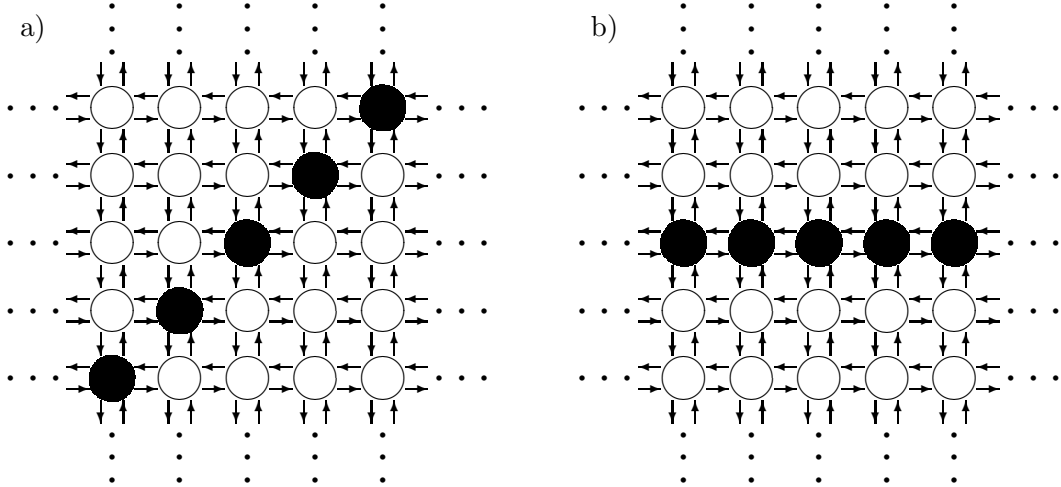}
\caption{Two examples of source-transitive graphs embedded in $\ZZ^2$ with an infinite motif ; a) a diagonal of sources; b) a horizontal array of sources.}
\label{fig:other-transitive}
\end{figure}
\end{center}

Finally, we extend naturally  our approach to a wide class of infinite  graphs that are called \emph{transitive}, but with a finite set of habitat qualities.  As phrased in \cite{TDW}, a graph is transitive if it globally ``looks the same'' from any vertex (including labels and edges). The vertices of the graph  are now labeled by  a countable set $D$ and $d_{PQ}$ gives the dispersal probability from patch $P \in D$ to patch $Q \in D$.  The type of patch $P$, denoted by  $i=\jmath(P) \in \{1, \ldots K\}$,  gives its quality  and the mean number of offspring per indivual in patch  $P$ is equal to $m_{\jmath(P)}$. In this paper, we will assume that the graph is only source-transitive, that is, transition probabilities and habitat qualities seen from any  patch of type  $1$ are identical.
 We specify the mathematical definition in Section \ref{transitivegraphs}.
 In particular, a practical example is given by  sources with the same quality  connected by corridors of identical  sinks and of the same length (see Figure \ref{fig:pipeline} for an example). A transitive  graph could also be  an (infinite) chessboard where whites are sinks and blacks are sources, the square lattice $\ZZ^2$ where sources have coordinates of type $(n,n)$ (diagonal) or of type $(n,0)$ (horizontal array), and so on (see Figures \ref{fig:transitive} and \ref{fig:other-transitive})...

\paragraph{Background.} In ecology and conservation biology, an important issue is to understand the effects of spatial heterogeneity on population persistence, in particular in the context of global climate change and wildlife preservation management. In particular, recent papers have focused on the factors (dispersal, temporal autocorrelation) that can allow for persistence even in sink metapopulations, see e.g. \cite{JansYosh, Gonz, RHB} and references therein.

In mathematics, the stochastic process describing the size of a density-independent, individual-based, stochastic population is called a \emph{branching process} \cite{asmu, AN}. Here, we consider multitype branching processes, where the type of an individual is given by its dwelling patch. When  the environment varies, the population can be  described by multitype branching processes in varying or random environment, see  \cite{athreyakarlin, kaplan, Tanny}  for mathematical references. 

The number of individuals located in  patch $i$  in generation $n$ is  denoted by $Z_n^{(i)}$. The process  $Z=(Z_n^{(i)}, \ i=1,\cdots,K, \   n \geq 0)$ is a \emph{multitype Galton--Watson process}. It is known from the mathematical literature \cite{asmu, AN} that either\footnote{we exclude in this paper  the  degenerate case when the size of population is constant with probability~1} the population becomes extinct or it grows exponentially. More specifically,  we see that $m_{i}d_{ij}$ is equal to the mean number of offspring of an individual living in patch $i$ which will land into patch $j$ in one time step, and therefore we call \emph{mean offspring matrix}  the matrix $A$ defined as
$$
A:=(m_{i}d_{ij} \  :  \ i,j=1,\cdots,K).
$$
The maximal eigenvalue (see e.g. \cite{Seneta}) of $A$ is called \emph{Malthusian growth rate}, or simply \emph{growth rate}. Indeed \cite{AN, asmu}, if $\rho\le 1$, the metapopulation dies out with probability $1$, and if $\rho>1$, the metapopulation can survive with positive probability, in which case 
$$
\frac{Z_n}{\rho^n}\stackrel{n\rightarrow \infty}{\longrightarrow} W ,
$$
where $W$ is  multidimensional finite and positive\footnote{under some integrability assumptions on the offspring numbers and the strong irreducibility (or the primitiveness) of   the matrix $A$. These conditions will be fulfilled here.}.

\paragraph{Method.}
We call random disperser a single individual who moves on the graph at discrete time steps following the dispersal probabilities. In other words, if $X_n$ denotes the position of such a random disperser after $n$ time steps, then $(X_n)$ is the Markov chain with transition probabilities
$$
\PP(X_{n+1}=j\mid X_n=i) = d_{ij}.
$$
The goal of this paper is to display persistence criteria, along with results regarding the asymptotic growth rate and the asymptotic fraction of time spent in each patch (by an individual taken at random in the surviving population) in terms of some characteristics of the motion of the random disperser. In contrast with the method involving the maximal eigenvalue of the mean offspring matrix, this one can yield quite simple, interpretable and partially explicit criteria. In addition, these criteria decouple the effect of dispersal and reproduction on population survival. This approach is still valid when the graph is an infinite transitive graph with finite motif, and partially valid even with infinite motifs.
In a number of remarks, we will also provide sufficient conditions for survival which are explicit. \\

\begin{table}[ht]%
\centering
\begin{tabular}{c|c}
 Notation & Interpretation  \\\hline
vertex $i$ & patch \\\hline
oriented edge  & probability of dispersal from patch $i$ to patch $j$ \\
 with weight $d_{ij}$ & \\ \hline
$m_i$ & mean number of offspring in a patch with habitat type $i$\\\hline
$M=m_1$ & mean growth rate in a source habitat (case $m_1>1\geq m_2=\cdots=m_K$)\\\hline
$m=m_2$ & mean growth rate in a sink habitat (case $m_1>1 \geq m_2=\cdots=m_K$)\\\hline
$p=d_{12}$ & probability of dispersal from the source to the sink (case of 2 patches) \\\hline
$q=d_{21}$ & probability of dispersal from the sink to the source (case of 2 patches)\\\hline
$X$ & random walk on the graph following the dispersal probabilities\\\hline
\end{tabular}
\caption{Notation.}
\end{table}

We assume that at least one loop-edge has a positive weight, that is, there is at least one patch in which the probability of staying put is nonzero. We further make the assumption that any patch can be visited with positive probability by the random disperser whatever its initial position. The previous two assumptions make the Markov chain $(X_n;n \ge 0)$ irreducible and aperiodic and the matrix $D$ (and $A$) primitive (or strongly irreducible). Actually, our results could be extended to this framework and even to the general case by lumping states together into classes but we will not develop this point here. 
We do not consider here the degenerate case when the  \emph{per capita} offspring number is constant and a.s. equal to $1$ in patch $1$. Finally, we assume that the offspring number $N_i$ of any individual living  in patch $i$ satisfies $\EE(N_i\log^+ N_i)<\infty$ (finite $N\log N$ moment). This guarantees the convergence of $Z_n/\rho^n$ to a non degenerate r.v. $W$ which is positive on the survival event.

\section{A first result on global persistence}

\subsection{General case}

We now give a criterion for metapopulation persistence in terms of the random disperser $X$. 
For that purpose, we assume from now on that the random disperser starts  in patch $1$  ($X_0=1$) and we denote by $T$ the first return time of the random disperser into patch $1$,  
$$
T:=\min\{n\ge 1: X_{n}=1\}.
$$
\begin{thm}
\label{eqn : persistence crit2}
The population persists with positive probability iff
$$
m_1\EE(\Pi_{i=1}^{T-1} m_{X_i}) > 1.
$$
\end{thm}
We observe that 
$$\EE(\Pi_{i=1}^{T-1} m_{X_i})=\Pi_{i=2}^K m_{i}^{S_{n}(i)},
$$
where $S_n(i):=\#\{ 1\leq k\leq T-1 : X_k=i\}$ is the \emph{time spent in habitat $i$ by time $T$}, that is, the number of times the random disperser has visited patches of habitat type $i$ until time  $T$.

\paragraph{Proof.}
Let $a$ (ancestor) be some  individual placed at time $0$ in  patch $1$. Define $Y_1$ as the number of offspring of $a$ staying put in patch $1$. Now for any integer $n\ge 2$, let $Y_n$ denote the number of descendants of $a$ at generation $n$ living in patch $1$ and whose ancestors  at generations $1, 2,\ldots, n-1$ have all lived \emph{outside} patch $1$. 
Then set
$$
Y:=\sum_{n\ge 1}Y_n,
$$
that can be seen as the total number of descendants of $a$ who live in patch $1$  for the first time in their lineage (except $a$).

In the theory of random trees, the set of individuals belonging to one of the $Y_n$ individuals for some $n$, is called a stopping line (see e.g. \cite{Chauvin}). By the branching property of our tree, the total numbers of descendants of each of the individuals of this stopping line who live in patch $1$  for the first time in their lineage, are independent and all follow the same law as $Y$. Therefore, the total number of descendants of $a$ living in patch $1$  is finite iff the branching process with offspring number distributed as $Y$ is finite, which is equivalent to the a.s. extinction of $Y$.  The bottomline is that the population persists with positive probability iff $\EE(Y)>1$.  Indeed, we have excluded the case when the per capita offspring is a.s. constant equal to one in patch $1$. \\

Let us compute $\EE(Y)$. We first note  that for every $i=1,\ldots,K$,  
$$\EE( Z_{n}^{(i)}  )=\sum_{j=1}^K \EE( Z_{n-1}^{(j)})m_{j}d_{ji},$$
where  $Z_{n}^{(i)}$ denotes the number of individuals located in patch $i$ at generation $n$. We prove easily by induction that the number of individuals $Y^{(i)}_n$ in patch $i$ at generation $n$ which have avoided patch $1$ at generations $k=1,2, \cdots, n-1$ satisfies
$$\EE(Y_n^{(i)})=\sum_{j=2}^K \EE(Y_{n-1}^{(j)})m_jd_{ji}=\sum_{j_1,\ldots, j_{n-1} \in \{2,\ldots K\}} d_{1j_1}d_{j_1j_2}\ldots d_{j_{n-2}j_{n-1}}d_{j_{n-1}i} m_1m_{j_1}\cdots m_{j_{n-2}}m_{j_{n-1}}.$$
As $Y_n=Y_n^{(1)}$, we get
\begin{eqnarray*}
\EE(Y_n)&=& \sum_{j_1,\ldots, j_{n-1} \in \{2,\ldots K\}} d_{1j_1}d_{j_1j_2}\ldots d_{j_{n-2}j_{n-1}}d_{j_{n-1}1} m_1m_{j_1}m_{j_2}\cdots m_{j_{n-2}}m_{j_{n-1}} \\
&=& m_1\EE(1_{ T=n}m_{X_1}m_{X_{2}}\cdots m_{X_{n-2}}m_{X_{n-1}}).
\end{eqnarray*}
Adding that  $Y=\sum_{n\geq 1} Y_n^{(1)}$, we have 
$$
\EE(Y)=   m_1\sum_{n \ge 1 } \EE(1_{ T=n}m_{X_1}\cdots m_{X_{n-1}}) =m_1\EE(m_{X_1}\cdots m_{X_{T-1}}).$$
This yields the result. \hfill$\Box$

\subsection{Case of two habitat types}
Let us focus now on the special case when there are $2$ habitat types and the source is  solely connected to sinks :
$$M:=m_1 >1, \qquad m:=m_2=\cdots=m_{K}<1.$$
We denote by 
$$p=\sum_{j=2}^K d_{1j}$$
 the probability of dispersing  for an individual living in patch $1$. The \emph{per capita} mean offspring number sent out from a source at each generation is $Mp$.
Let $\sigma$ be the time of \emph{first visit of a sink} by the random disperser
$$
\sigma:=\inf\{n\ge 0: X_n \ne 1\},
$$
so that $\sigma$ is a geometric random variable with success probability $p$.
Next, let $S$ denote the waiting time (after $\sigma$) before the random disperser visits a source (this source might or might \emph{not} be the initial source patch $X_0$)
$$
S:=\inf\{n\ge 0: X_{\sigma +n}=1\}.
$$
The duration $S$ can be seen as the \emph{time spent in sinks between two consecutive visits of sources}. In the two-vertex case, and in this case only, it is a geometric random variable with success probability $q$. The previous theorem reads as follows. 
\begin{prop}
\label{thm initial}
The population persists with positive probability iff
\begin{equation}
\label{eqn : persistence crit}
M(1-p) + eMp >1,
\end{equation}
where $e$ is the depleting rate due to the sink habitat in the graph, defined as
$$
e:=\EE(m^S)= \sum_{k\ge 1}m^k \PP(S=k).
$$
\end{prop}
To see that, we use  the first transition of the random disperser to get $\EE(\Pi_{i=1}^{T-1} m_{X_i})=(1-p)+p\EE(m^S)$.
\begin{rem}
If the average time spent in sinks has
$$
\EE(S)<\frac{M-1}{Mp(1-m)},
$$
then the population persists with positive probability. Indeed, the mapping $f:x\mapsto f(x)=\EE(x^S)$ is convex so 
$$
e=f(m)\ge 1+f'(1)(m-1)=1-(1-m)\EE(S)>1-\frac{M-1}{Mp}=\frac{1-M(1-p)}{Mp},
$$
which yields $eMp+M(1-p)>1$.
\end{rem}

\begin{rem}
Criterion \eqref{eqn : persistence crit} is equivalent to
\begin{equation}
\label{eqn : persistence crit 2}
\EE(M^\sigma m^S)>1.
\end{equation}
Both criteria are equivalent formulations of Theorem \ref{eqn : persistence crit2}.
\end{rem}

In contrast with the computation of a maximal eigenvalue, our criterion allows to decouple the roles of the mean offspring numbers $M$ and $m$ from that of the metapopulation structure itself (encompassed in $S$ and $p$, or $S$ and $\sigma$), that is, of reproduction and survival. 

One can check in some simple cases that criterion \eqref{eqn : persistence crit} or the condition that the maximal eigenvalue $\rho$ of $A$ exceeds unity are equivalent. For example, in the case when there are only one source and one sink $(K=2)$ and two vertices. Here the mean offspring matrix $A$ is
$$
A=
\left(
\begin{array}{cc}
M(1-p)& Mp\\
mq & m(1-q)	
\end{array}\right) .
$$ 
The characteristic polynomial $C$ of this square matrix is 
$$
C(x)=(M(1-p)-x)(m(1-q) -x)-Mmpq.
$$
Either $M(1-p)> 1$ and the population living in the source ensures the persistence. Or 
$M(1-p) \leq 1$ and the
quadratic polynomial is convex and has non negative  derivative at $1$. Thus, its leading eigenvalue is greater than 1 iff $C(1)<0$, which reads
$$
\frac{Mp}{1-M(1-p)}> \frac{1-m(1-q)}{mq} .
$$
We recover \eqref{eqn : persistence crit} since here $S$ is geometric with success probability $q$, which yields
$$
e=\sum_{k\ge 1} q(1-q)^{k-1}m^k=\frac{mq}{1-m(1-q)} .
$$
Notice that even in this simple case where $A$ is a $2\times 2$ matrix, the computation of the leading eigenvalue is cumbersome, and we have used a trick to explicitly specify the persistence criterion.

\subsection{Example with pipes of identical sinks}
Assume that the source  is a vertex of degree $2$ in the graph, connected to a \emph{left} sink and a  \emph{right} sink.  The probability of staying put on a source is still $1-p$, the probability of dispersing onto a left sink is $pL$, and the probability of dispersing onto a right sink is $pR$ (so that $L+R=1$). The sinks form a pipeline of $n$ adjacent sinks linking adjacent sources.  
The probability of staying put on a sink is always $s$,  the probability of dispersing from a sink onto one of its two neighboring sinks is $r$ in the left-to-right direction of the pipe, and $l$ in the right-to-left direction of the pipe (so that $q=l+r=1-s$). See Figure \ref{fig:cycle} for an example.\\
 \begin{center}
\begin{figure}[ht]
\unitlength 1.3mm 
\linethickness{0.4pt}
\begin{picture}(68,32.813)(-15,0)
\put(28,26.813){\circle{6}}
\put(37,26.813){\circle{6}}
\put(46,26.813){\circle{6}}
\put(28,17.813){\circle{6}}
\put(28,8.813){\circle{6}}
\put(46,17.813){\circle{6}}
\put(46,8.813){\circle{6}}
\put(34.8224,6.6354){\rule{4.3551\unitlength}{4.3551\unitlength}}
\multiput(35.7957,10.8781)(0,-4.893){2}{\rule{2.4086\unitlength}{.7628\unitlength}}
\multiput(36.3585,11.5284)(0,-5.714){2}{\rule{1.283\unitlength}{.2832\unitlength}}
\multiput(37.529,11.5284)(-1.4561,0){2}{\multiput(0,0)(0,-5.6425){2}{\rule{.3981\unitlength}{.2116\unitlength}}}
\multiput(38.0918,10.8781)(-2.8148,0){2}{\multiput(0,0)(0,-4.6157){2}{\rule{.6311\unitlength}{.4856\unitlength}}}
\multiput(38.0918,11.2512)(-2.5622,0){2}{\multiput(0,0)(0,-5.1402){2}{\rule{.3786\unitlength}{.2638\unitlength}}}
\multiput(38.0918,11.4025)(-2.4307,0){2}{\multiput(0,0)(0,-5.3577){2}{\rule{.2471\unitlength}{.1787\unitlength}}}
\multiput(38.358,11.2512)(-2.9565,0){2}{\multiput(0,0)(0,-5.0676){2}{\rule{.2406\unitlength}{.1912\unitlength}}}
\multiput(38.6105,10.8781)(-3.5699,0){2}{\multiput(0,0)(0,-4.4404){2}{\rule{.3489\unitlength}{.3102\unitlength}}}
\multiput(38.6105,11.0758)(-3.4538,0){2}{\multiput(0,0)(0,-4.7287){2}{\rule{.2329\unitlength}{.2031\unitlength}}}
\multiput(38.8469,10.8781)(-3.9179,0){2}{\multiput(0,0)(0,-4.3442){2}{\rule{.224\unitlength}{.214\unitlength}}}
\multiput(39.0651,7.6087)(-4.893,0){2}{\rule{.7628\unitlength}{2.4086\unitlength}}
\multiput(39.0651,9.9048)(-4.6157,0){2}{\multiput(0,0)(0,-2.8148){2}{\rule{.4856\unitlength}{.6311\unitlength}}}
\multiput(39.0651,10.4235)(-4.4404,0){2}{\multiput(0,0)(0,-3.5699){2}{\rule{.3102\unitlength}{.3489\unitlength}}}
\multiput(39.0651,10.6599)(-4.3442,0){2}{\multiput(0,0)(0,-3.9179){2}{\rule{.214\unitlength}{.224\unitlength}}}
\multiput(39.2628,10.4235)(-4.7287,0){2}{\multiput(0,0)(0,-3.4538){2}{\rule{.2031\unitlength}{.2329\unitlength}}}
\multiput(39.4382,9.9048)(-5.1402,0){2}{\multiput(0,0)(0,-2.5622){2}{\rule{.2638\unitlength}{.3786\unitlength}}}
\multiput(39.4382,10.171)(-5.0676,0){2}{\multiput(0,0)(0,-2.9565){2}{\rule{.1912\unitlength}{.2406\unitlength}}}
\multiput(39.5895,9.9048)(-5.3577,0){2}{\multiput(0,0)(0,-2.4307){2}{\rule{.1787\unitlength}{.2471\unitlength}}}
\multiput(39.7154,8.1715)(-5.714,0){2}{\rule{.2832\unitlength}{1.283\unitlength}}
\multiput(39.7154,9.342)(-5.6425,0){2}{\multiput(0,0)(0,-1.4561){2}{\rule{.2116\unitlength}{.3981\unitlength}}}
\put(40,8.813){\line(0,1){.3687}}
\put(39.977,9.182){\line(0,1){.1826}}
\put(39.949,9.364){\line(0,1){.1805}}
\put(39.909,9.545){\line(0,1){.1777}}
\put(39.859,9.722){\line(0,1){.1743}}
\put(39.797,9.897){\line(0,1){.1702}}
\put(39.725,10.067){\line(0,1){.1654}}
\put(39.643,10.232){\line(0,1){.16}}
\put(39.551,10.392){\line(0,1){.154}}
\put(39.449,10.546){\line(0,1){.1474}}
\multiput(39.337,10.694)(-.03005,.03508){4}{\line(0,1){.03508}}
\multiput(39.217,10.834)(-.03216,.03316){4}{\line(0,1){.03316}}
\multiput(39.088,10.967)(-.03414,.03112){4}{\line(-1,0){.03414}}
\put(38.952,11.091){\line(-1,0){.144}}
\put(38.808,11.207){\line(-1,0){.1508}}
\put(38.657,11.314){\line(-1,0){.1571}}
\put(38.5,11.411){\line(-1,0){.1628}}
\put(38.337,11.498){\line(-1,0){.1679}}
\put(38.169,11.576){\line(-1,0){.1723}}
\put(37.997,11.642){\line(-1,0){.1761}}
\put(37.821,11.698){\line(-1,0){.1792}}
\put(37.642,11.744){\line(-1,0){.1816}}
\put(37.46,11.777){\line(-1,0){.1834}}
\put(37.277,11.8){\line(-1,0){.1844}}
\put(37.092,11.812){\line(-1,0){.1848}}
\put(36.908,11.812){\line(-1,0){.1844}}
\put(36.723,11.8){\line(-1,0){.1834}}
\put(36.54,11.777){\line(-1,0){.1816}}
\put(36.358,11.744){\line(-1,0){.1792}}
\put(36.179,11.698){\line(-1,0){.1761}}
\put(36.003,11.642){\line(-1,0){.1723}}
\put(35.831,11.576){\line(-1,0){.1679}}
\put(35.663,11.498){\line(-1,0){.1628}}
\put(35.5,11.411){\line(-1,0){.1571}}
\put(35.343,11.314){\line(-1,0){.1508}}
\put(35.192,11.207){\line(-1,0){.144}}
\multiput(35.048,11.091)(-.03414,-.03112){4}{\line(-1,0){.03414}}
\multiput(34.912,10.967)(-.03216,-.03316){4}{\line(0,-1){.03316}}
\multiput(34.783,10.834)(-.03005,-.03508){4}{\line(0,-1){.03508}}
\put(34.663,10.694){\line(0,-1){.1474}}
\put(34.551,10.546){\line(0,-1){.154}}
\put(34.449,10.392){\line(0,-1){.16}}
\put(34.357,10.232){\line(0,-1){.1654}}
\put(34.275,10.067){\line(0,-1){.1702}}
\put(34.203,9.897){\line(0,-1){.1743}}
\put(34.141,9.722){\line(0,-1){.1777}}
\put(34.091,9.545){\line(0,-1){.1805}}
\put(34.051,9.364){\line(0,-1){.1826}}
\put(34.023,9.182){\line(0,-1){.9199}}
\put(34.051,8.262){\line(0,-1){.1805}}
\put(34.091,8.081){\line(0,-1){.1777}}
\put(34.141,7.904){\line(0,-1){.1743}}
\put(34.203,7.729){\line(0,-1){.1702}}
\put(34.275,7.559){\line(0,-1){.1654}}
\put(34.357,7.394){\line(0,-1){.16}}
\put(34.449,7.234){\line(0,-1){.154}}
\put(34.551,7.08){\line(0,-1){.1474}}
\multiput(34.663,6.932)(.03005,-.03508){4}{\line(0,-1){.03508}}
\multiput(34.783,6.792)(.03216,-.03316){4}{\line(0,-1){.03316}}
\multiput(34.912,6.659)(.03414,-.03112){4}{\line(1,0){.03414}}
\put(35.048,6.535){\line(1,0){.144}}
\put(35.192,6.419){\line(1,0){.1508}}
\put(35.343,6.312){\line(1,0){.1571}}
\put(35.5,6.215){\line(1,0){.1628}}
\put(35.663,6.128){\line(1,0){.1679}}
\put(35.831,6.05){\line(1,0){.1723}}
\put(36.003,5.984){\line(1,0){.1761}}
\put(36.179,5.928){\line(1,0){.1792}}
\put(36.358,5.882){\line(1,0){.1816}}
\put(36.54,5.849){\line(1,0){.1834}}
\put(36.723,5.826){\line(1,0){.1844}}
\put(36.908,5.814){\line(1,0){.1848}}
\put(37.092,5.814){\line(1,0){.1844}}
\put(37.277,5.826){\line(1,0){.1834}}
\put(37.46,5.849){\line(1,0){.1816}}
\put(37.642,5.882){\line(1,0){.1792}}
\put(37.821,5.928){\line(1,0){.1761}}
\put(37.997,5.984){\line(1,0){.1723}}
\put(38.169,6.05){\line(1,0){.1679}}
\put(38.337,6.128){\line(1,0){.1628}}
\put(38.5,6.215){\line(1,0){.1571}}
\put(38.657,6.312){\line(1,0){.1508}}
\put(38.808,6.419){\line(1,0){.144}}
\multiput(38.952,6.535)(.03414,.03112){4}{\line(1,0){.03414}}
\multiput(39.088,6.659)(.03216,.03316){4}{\line(0,1){.03316}}
\multiput(39.217,6.792)(.03005,.03508){4}{\line(0,1){.03508}}
\put(39.337,6.932){\line(0,1){.1474}}
\put(39.449,7.08){\line(0,1){.154}}
\put(39.551,7.234){\line(0,1){.16}}
\put(39.643,7.394){\line(0,1){.1654}}
\put(39.725,7.559){\line(0,1){.1702}}
\put(39.797,7.729){\line(0,1){.1743}}
\put(39.859,7.904){\line(0,1){.1777}}
\put(39.909,8.081){\line(0,1){.1805}}
\put(39.949,8.262){\line(0,1){.1826}}
\put(39.977,8.444){\line(0,1){.3687}}
\put(35.5,2.813){\line(1,0){3}}
\put(38.5,32.813){\line(-1,0){3}}
\put(52,16.313){\line(0,1){3}}
\put(52,7.313){\line(0,1){3}}
\put(52,25.313){\line(0,1){3}}
\put(22,19.313){\line(0,-1){3}}
\put(22,10.313){\line(0,-1){3}}
\put(22,28.313){\line(0,-1){3}}
\put(38.5,2.813){\vector(0,1){3}}
\put(35.5,32.813){\vector(0,-1){3}}
\put(52,19.313){\vector(-1,0){3}}
\put(52,10.313){\vector(-1,0){3}}
\put(52,28.313){\vector(-1,0){3}}
\put(22,16.313){\vector(1,0){3}}
\put(22,7.313){\vector(1,0){3}}
\put(22,25.313){\vector(1,0){3}}
\put(35.5,4.313){\line(0,-1){1.5}}
\put(38.5,31.313){\line(0,1){1.5}}
\put(50.5,16.313){\line(1,0){1.5}}
\put(50.5,7.313){\line(1,0){1.5}}
\put(50.5,25.313){\line(1,0){1.5}}
\put(23.5,19.313){\line(-1,0){1.5}}
\put(23.5,10.313){\line(-1,0){1.5}}
\put(23.5,28.313){\line(-1,0){1.5}}
\put(26.5,14.813){\vector(0,-1){3}}
\put(26.5,23.813){\vector(0,-1){3}}
\put(44.5,23.813){\vector(0,-1){3}}
\put(44.5,14.813){\vector(0,-1){3}}
\put(34,28.313){\vector(-1,0){3}}
\put(34,10.313){\vector(-1,0){3}}
\put(43,10.313){\vector(-1,0){3}}
\put(43,28.313){\vector(-1,0){3}}
\put(29.5,11.813){\vector(0,1){3}}
\put(29.5,20.813){\vector(0,1){3}}
\put(47.5,20.813){\vector(0,1){3}}
\put(47.5,11.813){\vector(0,1){3}}
\put(31,25.313){\vector(1,0){3}}
\put(31,7.313){\vector(1,0){3}}
\put(40,7.313){\vector(1,0){3}}
\put(40,25.313){\vector(1,0){3}}
\put(37,.75){\makebox(0,0)[cc]{\scriptsize $1-p$}}
\put(41,5.5){\makebox(0,0)[cc]{\scriptsize $pR$}}
\put(33,11.5){\makebox(0,0)[cc]{\scriptsize $pL$}}
\put(42,11.5){\makebox(0,0)[cc]{\scriptsize $l$}}
\put(32,5.5){\makebox(0,0)[cc]{\scriptsize $r$}}
\put(53,9){\makebox(0,0)[cc]{\scriptsize $s$}}
\put(21,9){\makebox(0,0)[cc]{\scriptsize $s$}}
\put(53,27){\makebox(0,0)[cc]{\scriptsize $s$}}
\put(42,29.5){\makebox(0,0)[cc]{\scriptsize $r$}}
\put(41,23.75){\makebox(0,0)[cc]{\scriptsize $l$}}
\put(49,22){\makebox(0,0)[cc]{\scriptsize $r$}}
\put(43,22.5){\makebox(0,0)[cc]{\scriptsize $l$}}
\end{picture}
\caption{A  pipeline where  $n$ identical sinks ($n=7$) connect the source to itself.}
\label{fig:cycle}
\end{figure}
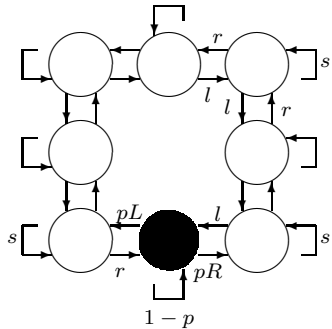
\end{center}

\begin{center}
\begin{figure}[ht]
\unitlength 1.3mm 
\linethickness{0.4pt}
\input{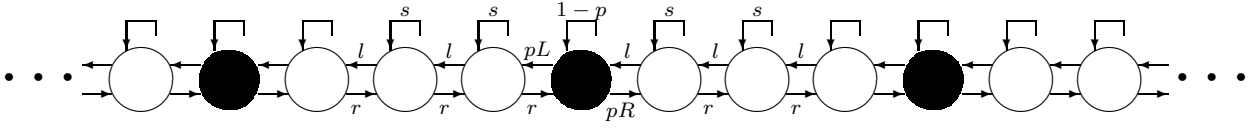}
\caption{A periodic pipeline where adjacent sources are separated by $n$ identical sinks ($n=3$).}
\label{fig:pipeline}
\end{figure}
\end{center}

This example will be directly  extended in the last Section to infinite transitive graphs, where pipelines of $n$ sinks periodically connect sources (see Figure \ref{fig:pipeline}). 

We can compute exactly the depleting rate $e$ of the above kind of source-transitive graph. Let $\lambda> 1>\mu$ be the two ordered solutions to
$$
mrx^2-(1-ms)x+ml=0\qquad x\ge 0.
$$
Then
$$
\lambda\mu = \frac{l}{r} \qquad\mbox{ and }\qquad \lambda+\mu = \frac{1-ms}{mr} .
$$
\begin{prop} The depleting rate $e$ is equal to
$$
e=\frac{\lambda^{n}-\mu^{n}}{\lambda^{n+1}-\mu^{n+1}}(L+R\lambda\mu)
+
\frac{\lambda-\mu}{\lambda^{n+1}-\mu^{n+1}}(R+L(\lambda\mu)^n) .
$$
\end{prop}
\begin{rem}
In the `chessboard' case ($n=1$), we recover 
$$
e=\frac{(1-s)m}{1-ms} .
$$
In the case when dispersal is isotropic $(l=r$), we get
$$
e=\frac{\lambda^{n}-\mu^{n}+\lambda-\mu}{\lambda^{n+1}-\mu^{n+1}} .
$$
Notice that in both previous cases, the depleting rate does not depend on $L$ or $R$.
In the case of one single source  and a large number of sinks ($n\tendinfty$), we get
$$
e=L\lambda^{-1} + R\mu.
$$
In the case of one single source and isotropic displacement, we then get $e=\lambda^{-1}=\mu$.
\end{rem}
\paragraph{Proof.} Consider a random walk $Y$ on $\{0,1,\ldots, n+1\}$, with displacement at each timestep being $-1$ with probability $l$, $0$ with probability $s$, and $+1$ with probability $r$. Let $T_i$ denote the first hitting time of $i$ by $Y$, and set $T:=\min(T_0, T_{n+1})$ as well as
$$
a_k:=\EE(m^{T}\mid Y_0=k).
$$
Returning to the random disperser $X$ on the graph, it is easily seen that 
$$
\EE(m^S\mid X_\sigma \mbox{ is a left neighbour})= a_1,
$$
while
$$
\EE(m^S\mid X_\sigma \mbox{ is a right neighbour})= a_n,
$$
so that
$$
e=La_1+Ra_{n}. 
$$
Computations of $a_1$ and $a_n$ rely on the following recurrence relationship
$$
a_k = ms a_k + ml a_{k-1}+mr a_{k+1} \qquad k\in\{1,\ldots, n\},
$$
with boundary conditions $a_0= a_{n+1}=1$. This relation is obtained easily by considering the first transition of the walk $X$. \hfill $\Box$

\section{Growth rate and habitat occupation frequencies}

\subsection{General case}

If $\mathbf u$ denotes an individual in generation $n$, we define $H_k(\mathbf u)$ as the patch  occupied by the ancestor of $\mathbf u$ in generation $k\le n$. Then, for every $i\in \{1,\ldots,K\}$,
$$
F_i(\mathbf u):=\frac{1}{n}\,\#\{ 0 \leq k \leq n : H_k(\mathbf{u})=i\}
$$
is the \emph{occupancy frequency of patch $i$ by the ancestral line of $u$}. We further denote by $\mathbf{U}_n$ an individual chosen randomly in the surviving population at generation $n$. We will see that the dispersal history of $\mathbf{U}_n$, as described by $(F_i(\mathbf{U}_n): i=1,\ldots,K)$ can be very different from that of a random disperser, since the ancestors of surviving individuals have better chance of having spent more time in sources than in sinks.

It is known that  the growth rate $\rho$ of the metapopulation is equal to  the maximal eigenvalue of $A$ (see \cite{asmu, AN,  KLPP}). Moreover, the asymptotic occupancy frequencies of $\mathbf{U}_n$ are deterministic and can be expressed as  the product  of the right and left maximal eigenvectors associated to $\rho$.  We refer to \cite{Jagers, JagersNermann, Hermissonal},  to Theorem 3.1 and 3.2 in  \cite{georgiibaake} in continuous time and to  \cite{KLPP} in discrete time.
In this section, we want to give an alternative characterization to these quantities in terms of the random disperser and show an application.
 
For that purpose, we use the transition matrix $D$ of the random disperser $X$ on the graph. We denote by $F_i(X_n)$ the \emph{occupancy frequency of patch $i$ by the random disperser $X$ by time $n$}
$$
F_i(X_n):=\frac{1}{n}\,\#\{ 0 \leq k \leq n : X_k=i\}.
$$
By assumption\footnote{the random walk $X$ on the graph is irreducible and aperiodic}, the random disperser has a a stochastic equilibrium on $V$, that we denote by $u=(u_i  : i=1,\ldots,K)$, which is the unique positive solution to $uD=u$. By the ergodic theorem, we also know that with probability 1,
$$
F_i(X_n)\stackrel{n\rightarrow \infty}{\longrightarrow} u_i .
$$
A typical single disperser will therefore occupy patch $i$ with asymptotic frequency $u_i$. This may not be the case of the ancestors of surviving individuals, whose paths must have favoured source patches. There is a trade-off between the pay-off in terms of fitness, gained by visiting source patches, and the cost in terms of likelihood, paid by deviating from the typical dispersal behaviour. This trade-off is particularly obvious if we consider the case of a perfectly sterile patch where the mean offspring is zero. In this case, the path followed by the ancestors of a surviving individual will necessarily have avoided this patch. Nevertheless, the asymptotic occupancy of this patch by a random disperser must be nonzero by assumption.

There is a way of quantifying both the cost and pay-off of deviating from the typical dispersal behaviour, that is, of having asymptotic occupancy frequencies $f=(f_i  : i=1,\ldots,K)$, where $f$ is a given element of the set $\cal F$ of non-negative frequencies on the graph
$$
\mathcal F:=\left\{f=(f_i :  i=1,\ldots,K) : \ f_i\geq 0, \ \sum_{i=1,\ldots,K} f_i=1\right\}.
$$
First, the probability that a random disperser has occupancy frequencies close to some given $f$ by time $n$ decreases exponentially with $n$ with rate $I(f)$, which can thus be interpreted as the cost of the $f$-occupancy scheme:
\begin{equation}
\label{def I}
I(f):=\sup\left\{\sum_{i=1,\ldots,K} f_i\log(v_i/(vD)_i) : v>>0\right\},
\end{equation}
where $v>>0$ denotes a positive vector, that is, $v_i>0$ for each $i=1,\ldots,K$.
Indeed, large deviations theory \cite{dembo, DenH} ensures that for any $\epsilon\ll 1$,  as $n\to\infty$, we have
\begin{eqnarray}
\label{frqc1}
\PP(   f_i-\epsilon \leq F_i(X_n) \leq f_i+\epsilon \mbox{ for all } i=1, \cdots, K) \asymp \exp(-n I(f)).
\end{eqnarray}
We refer the reader to forthcoming  Section \ref{proof} for a more  rigorous formulation. 
Taking $v=u$ in \eqref{def I} shows that $I(f)$ is of course always non-negative. This function is also convex.  When $f=u$, one can easily check that each partial derivative of $I$ is zero 
and  it can be proved that the supremum in \eqref{def I} is attained for $v=u$, so that $I(u)=0$. This was indeed expected, since $f=u$ is the natural occupancy scheme of the random disperser.

Second, the reproductive pay-off of $f$ can be defined as the fitness of a non-random disperser with given $f$-occupancy scheme, that is
$$
R(f):=\sum_{i=1,\ldots,K} f_i\log(m_{i}) .
$$
Indeed, the total size of a population of individuals  all adopting this dispersal behaviour can be seen to grow like $$\Pi_{i=1}^K m_i^{nf_i}=\exp(nR(f)).$$

Thus the cost (in terms of likelihood) for a population to follow some  occupancy scheme is quantified by $I$ and the reproductive pay off by $R$. The best strategy is to have an asymptotic occupancy frequency $\varphi$ which maximizes the difference $R-I$. If this optimum is positive, then the population  survives with positive probability. In addition, the ancestral line of a randomly chosen surviving individual will have visited patch $i$ with frequency $\varphi_i$. These results are stated below. The last assertion indicates that this optimal occupancy scheme $\varphi= (\varphi_i,i=1\ldots, K)$ is always different from the natural occupancy scheme $u$ of one single random disperser, except when all habitat types have the same quality.

\begin{thm} \label{rateped} The growth rate $\rho$ of the metapopulation is given by
\begin{eqnarray*}
\log (\rho) &=&\max\left\{ R(f) - I(f): f\in \mathcal F\right\}. 
\end{eqnarray*}
In addition, for any patch $i=1,\ldots,K$, conditional on the population being alive at time $n$, the occupancy frequency of patch $i$ by the ancestral line of a randomly chosen individual $\mathbf{U}_n$ in the surviving population at time $n$, converges to $\varphi_i$ in probability
$$ 
F_i(\mathbf{U}_n)\stackrel{n\rightarrow \infty}{\longrightarrow } \varphi_i,
$$
where the frequency vector $\varphi \in \mathcal F$  is uniquely  characterized by
$$
\log (\rho) = R(\varphi) - I(\varphi).
$$
The occupancy frequency $\varphi$ coincides with the stationary distribution $u$ of  $X$ (if and) only if
$$
m_1=m_2=\ldots =m_K.
$$
\end{thm}
In the same vein, we refer to Theorem 3.3 in \cite{georgiibaake} for a description of the lineage of surviving individuals for multitype branching processes in continuous time.

The proof is deferred to Section \ref{proof}. Strong irreducibility guaranteed by the assumptions given in the Introduction is useful (only) for the last part of the Theorem.

\subsection{Example with two patches} 

Let us apply this result for the simple case of two patches. Recall that $p$ is the probability to go from  patch $1$ to patch $2$ and $q$ the probability to go from patch $2$ to patch $1$ (see Figure
\ref{fig:two-patches}). 
We compute $I(f_1,f_2)$ which amounts to finding the maximum of 
$$u \in [0,1] \rightarrow f_1\log(u/(u(1-p-q)+q))+f_2\log((1-u)/(u(p+q-1)+1-q)).$$
We search the roots of its first derivative  : 
$$f_1\left( \frac{1}{u}-\frac{1-p-q}{(1-p-q)u+q}\right) +f_2\left(-\frac{1}{1-u}+\frac{p+q-1}{(p+q-1)u+1-q}\right)=0.$$
Now we focus on the case  $p+q=1$, assumption which is made in \cite{JansYosh} and makes computations easier (the general case requires to solve a third degree polynomial equation).  Roughly speaking, this corresponds to letting a fraction $q$ of the population live in patch $1$ (and  $p=1-q$ live in patch 2).
Then we check directly from the last displayed equation that
$$u=f_1, \quad 1-u=f_2=1-f_1$$ 
and $$I(f_1,f_2)=f_1\log(f_1/q)+f_2\log(f_2/(1-q)).$$
To compute $\max\left\{ R(f) - I(f): f\in \mathcal F\right\}$, we just need to maximize
$$f_1 \in [0,1]\rightarrow f_1\log(m_1)+(1-f_1)\log(m_2)-f_1\log(f_1/q)-(1-f_1)\log((1-f_1)/(1-q)). $$
 We search again the roots of the first derivative
$$\log(m_1/m_2)-\log(f_1/q)-1+\log((1-f_1)/(1-q))+1=0,$$
which yields the habitat occupation frequency of patch $1$ : 
$$ \frac{1-\varphi_1}{\varphi_1}=\frac{1}{\varphi_1}-1=\frac{(1-q)m_2 }{qm_1}.  $$
We get 
$$\varphi_1 =\frac{qm_1 }{qm_1+(1-q)m_2}=1-\varphi_2$$
and also recover the expected result $\rho=qm_1+(1-q)m_2$.

\section{Fluctuating  environments}

\subsection{General setting}

We now enrich our model with a fluctuating environment. The environment is embodied by a certain value $w$ which belongs to a finite set of states. We assume that the environment affects simultaneously all patches, but not necessarily  in the same way. We keep on assuming a simple asexual life cycle with discrete non-overlapping generations and no density-dependence. Now the environment is assumed to affect the reproduction scheme, but not the dispersal behaviour. Specifically, at each time step, conditional on the state $w$ of the environment, individuals reproduce independently according to some distribution which depends on the habitat type of their dwelling patch. We denote by $m_i(w)$ the mean offspring number of individuals dwelling in patch $i$ when the environment is in state $w$. 

Except in the last subsection, we will assume that the environment alternates periodically at each time step between two states. Actually, the same method would allow to deal with any finite number of environmental states varying periodically, which allows for the modeling of seasonal effects.  In contrast with \cite{JansYosh} who deal with large populations sending out at each generation given fractions of their offspring to each patch (parent-independent migration, also known as island model), we consider discrete populations with stochastically varying size, individually stochastic dispersal behaviours, and allow for metapopulation structure of the stepping-stone type. In this more general setting, we come to the same conclusion that in fluctuating environments, the population may survive in sinks only. This means that in the absence of dispersal, the population in each patch would become extinct. But positive dispersal  probabilities may allow for global survival. 

We call $e_1$ and $e_2$ the two possible states of the environment, so now we have $2K$  habitat qualities $m_i(e_j)$, for $i=1,\ldots, K$ and $j=1,2$. Now  the process $Z=(Z_n^{(i)},   i=1,\ldots,K ;n\ge 0)$ is no longer a homogeneous Markov chain but  a multitype branching process in varying environment. However, restricting the observation of the metapopulation to times when the environment is in the same state allows to adapt the arguments of the previous section. Indeed, $(Z_{2n}^{(i)}, i=1,\ldots,K;n\ge 0)$ is  a multitype branching process with mean offspring matrix $A$ with generic element
$$
a_{ij}=\sum_{k=1,\ldots,K} m_{i}(e_1)d_{ik}m_{k}(e_2)d_{kj}  \qquad i,j=1,\ldots,K.
$$
In the following subsection, we treat the case of two patches and determine the global persistence criterion. We then handle the general case using the random disperser $X$.

\subsection{Example with two patches and two periodic environments}

For simplicity, even if the environment is now variable, the two patches are still called respectively the source (patch 1) and the sink (patch 2). The mean number of offspring in the source are denoted by $M_1=m_1(e_1)$ and $M_2=m_1(e_2)$. In the sink, they are denoted by $m_1= m_2(e_1)$ and $m_2=m_2(e_2)$.

\begin{thm}
A necessary and sufficient condition for global persistence is
$$
M_1M_2(1-p)^2+(M_1m_2+m_1M_2)pq+m_1m_2(1-q)^2 > \min\big(2,1+ M_1M_2m_1m_2(1-p-q)^2\big).
$$
\end{thm}
\begin{rem}
It is easy to find examples where both patches are sinks when forbidding dispersal but the metapopulation survives with positive probability in the presence of dispersal. Indeed, in the absence of dispersal, each patch is a sink if (and only if) $M_1M_2\leq 1$ and $m_1m_2\leq 1$. Assuming for example that $p=q=1/2$ and $m_1=m_2=m$, the global survival  criterion becomes $M_1M_2+m(M_1+M_2)+m^2 > 4$, which holds as soon
 as $m(M_1+M_2) > 4$. 
\end{rem}

\paragraph{Proof.} The mean offspring matrix of $(Z_{2n};n\ge 0)$ is given by 
\begin{eqnarray*}
a(1,1)&=& M_1M_2(1-p)^2+M_1m_2pq \\
a(1,2)&=& M_1m_2p(1-q)+M_1M_2(1-p)p \\
a(2,1)&=&m_1M_2q(1-p)+m_1m_2(1-q)q \\
a(2,2)&=&m_1m_2(1-q)^2+m_1M_2qp.
\end{eqnarray*}
The maximum eigenvector of the matrix $A=(a(i,j) : 1\leq i,j\leq 2)$ is the largest root
 of the polynomial
 $$
 x^2-(a(1,1)+a(2,2))x+a(1,1)a(2,2)-a(1,2)a(2,1).
 $$
So it is less than $1$ iff
$$
a(1,1)+a(2,2)+\sqrt{(a(1,1)-a(2,2))^2+4a(1,2)a(2,1)}\leq 2
$$
Then the criterion for a.s. extinction of $(Z_{2n} : n \in \NN)$ is 
$$
a(1,1)+a(2,2)\leq 2 \quad \text{and} \quad (a(1,1)-a(2,2))^2+4a(1,2)a(2,1)\leq (2-a(1,1)-a(2,2))^2.
$$
The second inequality becomes
$a(1,1)+a(2,2)\leq 1+a(1,1)a(2,2)-a(1,2)a(2,1)$, which gives
$$
M_1M_2(1-p)^2+M_1m_2pq+m_1m_2(1-q)^2+m_1M_2qp\leq 1+ M_1M_2m_1m_2(1-p-q)^2.
$$
This completes the proof.\hfill $\Box$

\subsection{Global persistence for more than two patches}
Here, we extend the previous result to the case of a general, finite graph. We want to  state a global survival criterion which generalizes Theorem  \ref{eqn : persistence crit2} to periodic environments. Assume again that the random disperser starts at time $0$ in patch $1$  and set $T$ the first \emph{even} time when the random disperser goes back to habitat $1$
$$
T:=\min\{n\ge 1: X_{n}=1 \mbox{ and $n$ is even}\}.
$$

\begin{thm}
The population persists with positive probability iff
$$
m_1\EE(\Pi_{i=1}^{T-1} m_{X_i}(w_i)) > 1,
$$
where the sequence $(w_i : i\geq 1)$ can take one of the two values $(e_1,e_2,e_1,\ldots)$ or $(e_2,e_1,e_2,\ldots)$, depending on the initial environment. 
\end{thm}
This theorem can be proved by a direct adaptation of the proof of Theorem \ref{eqn : persistence crit2} replacing  $Z_n$ with $Z_{2n}$.  

\subsection{Rate of growth and habitat occupation frequency }

The generalization to periodic environments of the results of the previous section can be achieved by changing the state-space $\{1,\ldots,K\}$ of the random disperser to the state-space of oriented edges of the graph, i.e., ordered pairs of vertices 
$$
 \mathcal{E}:= \{1,\ldots,K\}^2.
 $$ 
Denote by $B$ the transition matrix of the Markov chain $(X_{2n}, X_{2n+1};n\ge 0)$, which indeed takes values in $ \mathcal{E}$. Then denote by $\mathcal{F}$ the set of frequencies indexed by $\mathcal{E}$
$$
\mathcal{F}:=\big\{(f_E, E\in \mathcal{E}) : \ f_E\geq 0, \quad  \sum_{E\in \mathcal{E}} f_E=1\big\},
$$
and define the new cost function $I:\mathcal{F}\rightarrow \RR$ as
$$
I(f):= \sup\left\{\sum_{E\in \mathcal{E}} f_E\log(v_E/(vB)_E) : v>>0\right\},
$$
where $v$ denotes a non-negative vector indexed by $\mathcal{E}$, such that $v>>0$, that is, $v_E>0$ for all $E\in \mathcal{E}$. Also define the new pay-off function  $R:\mathcal{F}\rightarrow \RR$ as
$$
R(f):=\sum_{E=(i,j) \in \mathcal{E}} f_E\log(m_{i}(e_1)m_{j}(e_2)).
$$
We can also provide an expression of  $I$ in terms of the entropy  function using Theorem 3.1.13 in \cite{dembo}.
The generalization of Theorem \ref{rateped} can be stated as follows.
\begin{thm} \label{ratetwo} The growth rate $\rho$ of the metapopulation is given by
\begin{eqnarray*}
2\log (\rho) &=&\max\{ R(f) -I(f) : f\in \mathcal F\}
\end{eqnarray*}
In addition, for any patch $i\in \{1,\ldots,K\}$, conditional on the population being alive at time $n$, the frequencies of occupation of patch $i$ by the ancestral line of a randomly chosen individual $\mathbf{U}_n$ in the surviving population at time $n$, converges  in probability :
$$
F_j(\mathbf{U}_n)\stackrel{n\rightarrow \infty}{\longrightarrow } \sum_{i \in \{1,\ldots,K\}} \varphi_{i,j}, 
$$
where the vector $(\varphi_{i,j} : (i,j)\in\mathcal{E})$ is characterized by
$$\log (\rho) =  R(\varphi) -I(\varphi).$$
\end{thm}
 The proof follows that of Theorem  \ref{rateped}, with now 
$$\EE(\vert Z_{2n+1} \vert )= \EE\left(\Pi_{i=0}^{n} m_{X_{2i}} m_{X_{2i+1}}\right)=\EE\left(\Pi_{i=1}^K\Pi_{j=1}^K m_i(e_1)^{S_n^{(1)}(i)}m_j(e_2)^{S_n^{(2)}(j)} \right),
$$
where $\vert Z_{n} \vert$ is the total number of individuals in source patches at generation $n$ and
$$
S_n^{(1)}(i)=\#\{ k\leq n : X_{2k}=i\}
,\quad S_n^{(2)}(i)=\#\{ k\leq n : X_{2k+1}=i\}
.
$$

\paragraph{Example with two patches.}

We  extend the computations of the previous section to periodic environments. We focus on the case $p+q=1$, which ensures that the mean fraction of individuals in patch $1$ is equal to $q$. We recall that this is an assumption made in \cite{JansYosh}.    \\
Let us first compute $I(f):= \sup\left\{\sum_{E\in \mathcal{E}} f_E\log(v_E/(vB)_E) : v>>0\right\}$ where
$$B_{(i,j)(k,l)}=d_{jk}d_{kl}=\epsilon_k\epsilon_l, \qquad i,j,k,l =1,2,$$
where $\epsilon_1=q=1-p, \epsilon_2=p=1-q$. Then, assuming that $\sum_{E\in \mathcal{E}} v_E=1$,
$$vB_{k,l}=\sum_{E\in \mathcal{E}} v_E \epsilon_k\epsilon_l=\epsilon_k\epsilon_l,$$
so that
$$I(f)= \sup\left\{\sum_{E\in \mathcal{E}} f_E\log\left(\frac{v_E}{\epsilon_{E_1}\epsilon_{E_2}}\right) : v>>0,  \sum_{E\in \mathcal{E}} v_E=1 \right\}$$
and
$$I(f)=\sum_{E\in \mathcal{E}} f_E\log\left(\frac{f_E}{\epsilon_{E_1}\epsilon_{E_2}}\right).$$
Then $R(f)-I(f)=\sum_{E\in \mathcal{E}} f_E\log\left(m_{E_1}(e_1)m_{E_2}(e_2)\epsilon_{E_1}\epsilon_{E_2}/f_E\right)$. Writing $f(2,2)=1-\sum_{E \in \mathcal{E}-(2,2)} f_E$ and using that  the partial derivatives of $R-I$ equal $0$ in $\varphi$,  we have :
$$ \log\left(m_{E_1}(e_1)m_{E_2}(e_2)\epsilon_{E_1}\epsilon_{E_2}\right)+1-\log(\varphi_E)-\left[\log\left(m_{2}(e_1)m_{2}(e_2)p^2)\right)+1-\log(\varphi(2,2))\right]=0$$
which gives the  habitat occupation frequencies
$$\varphi_E=\epsilon_{E_1}\epsilon_{E_2}m_{E_1}(e_1)m_{E_2}(e_2)\left[\sum_{E \in \mathcal{E}} \epsilon_{E_1}\epsilon_{E_2}m_{E_1}(e_1)m_{E_2}(e_2)\right]^{-1}.$$
We can also now deduce the growth rate $\rho$. Using that 
$$\sum_{E \in \mathcal{E}} \epsilon_{E_1}\epsilon_{E_2}m_{E_1}(e_1)m_{E_2}(e_2)=q^2M_1M_2+qpM_1m_2+pqm_1M_2+p^2m_1m_2, \qquad \sum_{E \in \mathcal{E}} f_E=1,$$ we get 
$$ R(\varphi)-I(\varphi)=\frac{1}{2}\log(q^2M_1M_2+qpM_1m_2+pqm_1M_2+p^2m_1m_2)$$
and
$$\rho=(qM_1+pm_1)^{1/2}(qM_2+pm_2)^{1/2}.$$
This is the same growth rate as the one computed in \cite{JansYosh}.


\subsection{Some comments on random environments}
A more natural way of modeling fluctuating environment in ecology is to assume random rather than periodic environment. The approach developed for periodic environments cannot be extended to random environments directly. Indeed, since the environment affects the whole metapopulation simultaneously, the randomness of environments correlates reproduction success in different patches. The process $(Z_n^{(i)}, i=1,\ldots,K;n\ge 0)$ counting the population sizes on each patch is now a \emph{multitype branching process in random environment} (MBPRE).


Let us denote by $A(w)$ the mean offspring matrix (involving dispersal) in environment $w$. Specifically, the generic element $a_{ij}(w)$ of $A(w)$ is the mean offspring number of a typical individual dwelling in patch $i$ sent out to patch $j$ by dispersal, when the environment is $w$, so that
$$
A_{ij}(w)=m_{i}(w)\,d_{ij}.
$$
We will now assume that the state-space of environments is finite and that the sequence $(w_n: n \geq 0)$ of environment states through time is a \emph{stationary, ergodic sequence}, possibly autocorrelated, in the sense that the states need not be independent. 
Under this assumption, it is proved \footnote{under the further assumption $\EE(\log^+\parallel A(w_0)\parallel )<\infty$, where expectation is taken w.r.t. the environment} in \cite{FK} that the limit $\gamma$ of the sequence
$$ 
\frac{1}{n}\, \log\parallel  A(w_{n})A(w_{n-1})\ldots A(w_0)\parallel 
$$
exists with probability 1 and is deterministic, where $\parallel.\parallel$ denotes the maximum row sum of the matrix. This is interesting to us because it is further shown\footnote{again under some further assumptions,  see Annex.} in \cite{athreyakarlin, kaplan, Tanny}, that the extinction criterion and the growth rate of this MBPRE are respectively given by the sign and the value of $\gamma$, more specifically, $\gamma=\log \rho$. 

Unfortunately, this does not give a very explicit condition for global persistence. But again using the random disperser, we can give some sufficient conditions for survival. For simplicity, we turn our attention to the example of two patches and two environments $e_1$ and $e_2$. At any time step, the probability that the environment is in state $e_1$ is denoted by $\nu \in (0,1)$ (so the probability that the environment is in state $e_2$ is $1-\nu$). We show again that  the population may survive in sinks only.
As in the case of periodic environments, the mean number of offspring in the first patch (the source) is denoted by  $M_1=m_1(e_1)$ and $M_2=m_2(e_1)$. In the second patch (sink), they are denoted by $m_1= m_2(e_1)$ and $m_2=m_2(e_2)$. 


For the sake of simplicity, we state the results for two special cases of ergodic environments: the case when the states are independent at each time step (with probability $\nu$ and $1-\nu$), and the case when the sequence is a Markov chain. We denote by $\alpha$ the transition from $e_1$ to $e_2$ and by $\beta$ the transition from $e_2$ to $e_1$. Then it is well-known that $\nu=\beta/(\alpha+\beta)$. The case of independent environments is recovered when $\alpha+\beta=1$. Note that as soon as $\alpha+\beta\not=1$, the sequence of environment states is auto-correlated.

\begin{prop} We have the following lower bound for the growth rate of the metapopulation.
$$
\log(\rho)\geq \nu\log(M_1)+(1-\nu) \log(m_{2})+ \nu\alpha\log(pq)+\nu(1-\alpha)\log(1-p)+(1-\nu)(1-\beta)\log(1-q).
$$
\end{prop}
\begin{rem}
Observe that this lower bound does not depend on $M_2$ and $m_1$. Again one can display examples where both patches are sinks when forbidding dispersal but the metapopulation survives with positive probability in the presence of dispersal. In the absence of dispersal, each patch is a sink if (and only if) $M_{1}^{\nu}M_{2}^{1-\nu}\le 1$ and $m_{1}^{\nu}m_{2}^{1-\nu}<1$. Actually one can manage to keep $\gamma >0$ while  $M_{1}^{\nu}M_{2}^{1-\nu}<1$, $m_{1}<1$ and $m_{2}<1$, for example with $M_{2}$ small and $M_1$ large for some fixed $p,q, m_1,m_2$. This corresponds to $e_2$ being a catastrophic environment in the source patch but the population survives  in patch $2$ when a  catastrophe occurs. 
\end{rem}


\paragraph{Proof.} We consider only  the subpopulation avoiding patch $1$ when the environment is equal to $e_2$. This means that this population reproduces with mean offspring number $M_1(1-p)$ in patch $1$ while the environment is  $e_1$. Each time the environment $e_2$ occurs, we consider the part of this population which has dispersed to patch $2$. This corresponds to a mean offspring number of $M_1p$. This population then stays in patch $2$ and  reproduces with mean offspring number $m_2(1-q)$ until the environment is again equal to $e_1$. We then consider the part of this population which goes back to patch $1$. This corresponds to a mean offspring number of $m_2q$.

Thus the patch of the ancestors of the individuals we keep is equal to $1$ (resp. $2$) if it lived  in environment $e_1$ (resp. $e_2$). Then at time $n$, the mean size of the population we consider is  equal to 
$$M_1^{N_1(n)}m_{2}^{N_2(n)}(1-p)^{N_{11}(n)}(1-q)^{N_{22}(n)}p^{N_{12}(n)}q^{N_{21}(n)}$$
where $N_i(n)$ ($i\in\{1,2\}$) is the number of times before generation $n$ when the environment is equal to $e_i$ 
 and $N_{ij}(n)$ ($i,j\in\{1,2\}$) is the number of one-step transitions of the environment from  $e_i$ to type $e_j$ until time $n$. By ergodicity, we know that  these quantities have  deterministic frequencies asymptotically. In the case of a Markovian sequence of environments, as $n\rightarrow \infty$,
$$
N_1(n)\sim \nu n, \quad N_2(n)\sim (1-\nu) n,
$$
and 
$$ 
N_{11}(n)\sim \nu(1-\alpha) n, \quad N_{12}(n)\sim \nu\alpha n,\quad N_{21}(n)\sim (1-\nu)\beta n, \quad N_{22}(n)=(1-\nu)(1-\beta)n.
$$
Using the growth of this  particular part of the whole population directly gives us a lower bound for $\gamma$:
$$
\gamma\geq \nu\log M_1+(1-\nu) \log m_{2}+ \nu\alpha\log p+(1-\nu)\beta\log q+\nu(1-\alpha)\log(1-p)+(1-\nu)(1-\beta)\log(1-q).
$$
Noticing that $(1-\nu)\beta=\nu \alpha$ completes the proof.  \hfill$\Box$

\begin{rem}
We could improve these results by considering more sophisticated strategies. For example, we could consider the subpopulation which stays in patch $1$ if (and only if) the number of consecutive catastrophes is less than $k$  and then optimize over $k$.\\
Actually, the proof relies on a stochastic coupling. Roughly speaking, the subpopulation we consider avoids the bad patches at the bad times and follows a (one type) branching process in random environment $e_{11}, e_{12}, e_{21}$ and $e_{22}$ respectively with stationary probabilities $\nu(1-\alpha)$, $\nu\alpha$, $(1-\nu)\beta$ and $(1-\nu)(1-\beta)$and mean offspring $M_1(1-p)$, $M_1p$, $m_2q$ and $m_2(1-q)$.  \\
Observe also that we can derive  a lower bound using the permanent of the mean matrix $M$ of the MBPRE from Proposition $2$ in \cite{BS}. But this  lower   bound is not relevant for understanding the survival event in sinks only.
\end{rem}

\section{Metapopulation on infinite transitive graphs}
\label{transitivegraphs}

We now turn our attention to infinite graphs with  a special property called \emph{transitivity}. We label the patch of the graph with a countable set $D$. 
Each patch $P\in D$ has a type $i=\jmath(P) \in \mathbb N$  which gives  its quality. That is, the mean number of offspring in patch $P$ is equal to $m_{\jmath(P)}$.
The transitivity assumption means here  that the graph, including  dispersal probabilities and habitat types, looks the same from any patch of type $1$.

\begin{center}
\begin{figure}[h]
\unitlength 1mm 
\linethickness{0.4pt}
\input{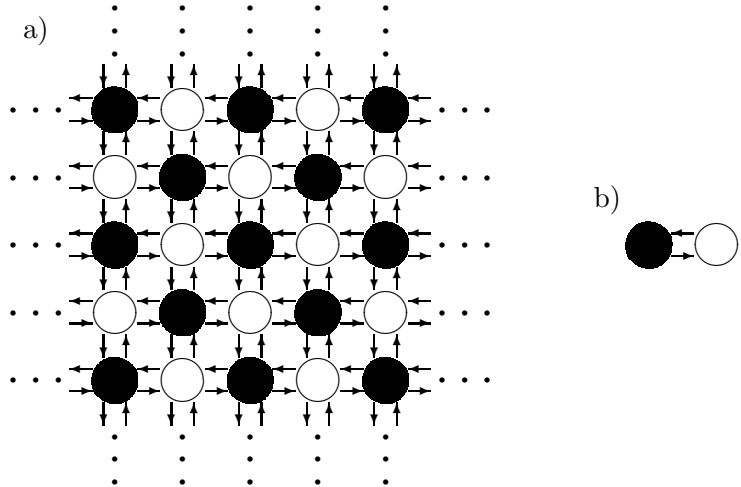}
\caption{a) The chessboard (with periodic arrow labels - not represented) is a source-transitive graph that can be collapsed into: b) a two-vertex graph (loop edges are not represented).}
\label{fig:chessboard}
\end{figure}
\end{center}

\begin{center}
\begin{figure}[h]
\unitlength 1mm 
\linethickness{0.4pt}
\input{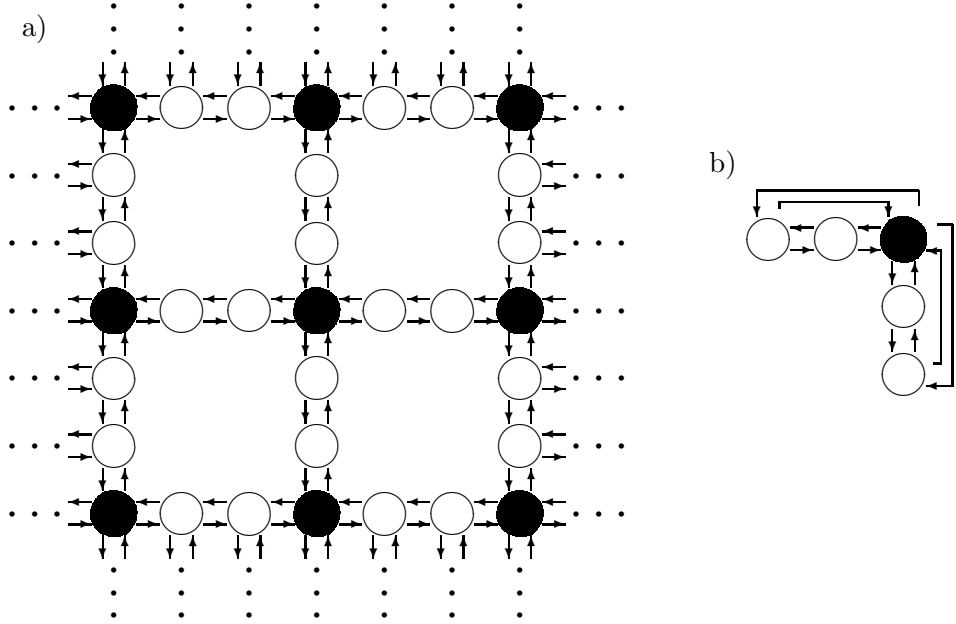}
\caption{a) Star sources with $2d$ pipelines of sinks (and periodic arrow labels - not represented) form a source-transitive graph (here $d=2$) that can be collapsed into: b) a single source with $d$ pipelines  (loop edges are not represented).}
\label{fig:stars}
\end{figure}
\end{center}

Let us provide the reader with some examples. An infinite linear periodic array of patches is a transitive graph (see Figure \ref{fig:array} for example). This means that two consecutive sources are separated by the same sequence of  habitat qualities. The patches can then be labeled by $i\in D=\ZZ$ and the type of patch $i$ is equal to $f(i)$, for some integers $N,K>0$ and a function  $f : \ZZ \rightarrow [1,K]$ such that $f(i+N)=f(i)$ for every $i\in \ZZ$. The chessboard is a transitive graph (see Figure \ref{fig:chessboard}). Star sources with $2d$ pipelines of sinks form a transitive graph (see Figure \ref{fig:stars} for $d=2$).

Let us now specify mathematically these definitions. We consider here graphs whose vertices are patches which can be  labeled by $P \in D$, the oriented edges from $P$ to $Q$  are weighed by $d_{PQ}$  and 
the type of vertex $P$ is denoted by  $\jmath (P)$. A  mapping $T$ of the graph is called  an isomorphism if it conserves the types of the vertices as well as the weights of the oriented edges: $T$ is a bijection of $\NN$ such that  for all $P,Q \in \NN$,
$$ d_{T(P)T(Q)}=d_{PQ}, \qquad \jmath (T(P))=\jmath (P).$$
The associated equivalence relation $\sim$ between the patches of the graph is defined by 
$$P \sim P' \qquad \text{iff there exists an isomorphism } T \text{ of the graph such that } T(P)= P'.
$$
The class of a patch $P$ is defined as $Cl(P)=\{P' : P'\sim P\}$. Every patch of this class has the  type of  $P$. 
The graph is then called \emph{source-transitive} if  all patches of type $1$ belong to a single class.
The collection of the distinct classes $(Cl(i) : i \in V)$ form a partition of the patches of the graph. Such  subsets $V$ of patches are called  \emph{motifs}.   With a slight abuse of notation, the transition probabilities on  a motif $V$ are denoted by $(d_{PQ} : P\in V, Q\in V)$ where
$$d_{PQ}= \sum_{Q'\in Cl(Q)} d_{PQ'}$$
which does not depend on the element $P \in Cl(P)$. 

The initial graph can be seen as a family of copies of a motif properly connected. Observe that not all transitive graphs have a finite motif. For example, the cases illustrated by Figure \ref{fig:transitive} can be collapsed into a finite motif, but those given in Figure \ref{fig:other-transitive} (a diagonal or a horizontal array of sources in $\ZZ^2$) cannot. Indeed, Figure \ref{fig:chessboard} shows that the infinite chessboard in $\ZZ^2$ can be collapsed into a $2$-vertex graph and Figure \ref{fig:stars} shows that star sources with $2d$ pipelines can be collapsed into a finite motif with one source and $d$ pipelines of sinks.
\begin{table}[ht]%
\centering
\begin{tabular}{c|c}
 Notation & Interpretation  \\\hline
vertex $P$ & patch \\\hline
set $V$ of vertices & motif \\\hline
type of $P$ $\jmath(P)=i$ & habitat quality of patch $P$ is $i$ \\\hline
oriented edge with weight $d_{PQ}$ & probability of dispersal from patch $P$ to patch $Q$ \\\hline
$m_i$ & mean number of offspring in a patch with habitat type $i$\\\hline
\end{tabular}
\caption{Notation.}
\end{table}
We can easily extend the results of the previous section to infinite transitive graphs by considering the random disperser $X_n$ on the graph which follows the dispersal probabilities $(d_{PQ} : P,Q \in D^2)$.
In that purpose, we assume that the random disperser starts  in a patch of type $1$  ($\jmath(X_0)=1$) and we denote by $T$ the first return time of the random disperser into a (possibly different) patch of type $1$. 
$$
T:=\min\{n\ge 1: \jmath(X_{n})=1\}.
$$
\begin{thm}
The population persists with positive probability iff
$$
m_1\EE(\Pi_{i=1}^{T-1} m_{\jmath(X_i)}) > 1.
$$
\end{thm}
The number of individuals located in   patch $P$  in generation $n$ is  denoted by $N_n^{(P)}$. Collapsing the graph into some motif $V$, we denote by 
$$Z_n^{(P)}= \sum_{ P' \in Cl(P)} N_n^{(P)}$$
the total number of individuals in some habitat of the class of patch $P\in V$. 
In the case when the motif is finite, $Z=(Z_n^{(P)}, \ P\in V;  n \geq 0)$ is a multitype Galton--Watson process.
We see that $m_{\jmath(P)}d_{PQ}$ is equal to the mean number of offspring of an individual living in patch $P$ which will land into patch $Q$ in one time step, and therefore we call mean offspring matrix  
$$
A:=(m_{\jmath(P)}d_{PQ} \  :  \ P,Q\in V),
$$
Then, everything  happens as if the metapopulation evolves on a finite graph corresponding to the motif. We then derive the following counterpart of Theorem \ref{rateped}. We denote now by $D$ the transition matrix of the Markov chain $\jmath(X_n)$, and still use for $f=(f_P : P\in V)$
$$
I(f):=\sup\big\{\sum_{P\in V} f_P\log(v_P/(vD)_P) : v>>0\big\},
$$
and 
$$
R(f):=\sum_{P \in V} f_P\log(m_{\jmath(P)}) .
$$
\begin{thm} We assume that the graph is transitive and has a finite motif. \\
 The growth rate $\rho$ of the metapopulation is given by
\begin{eqnarray*}
\log (\rho) &=&\max\left\{ R(f) - I(f): \sum_{P\in V} f_P=1, \ f_P \geq 0\right\}. 
\end{eqnarray*}
In addition, for any $P\in V$, conditional on the population being alive at time $n$, the occupancy frequency of a patch with type $\jmath(P)$ by the ancestral line of a randomly chosen individual $\mathbf{U}_n$ in the surviving population at time $n$, converges to $\varphi_i$ in probability
$$ 
F_i(\mathbf{U}_n)\stackrel{n\rightarrow \infty}{\longrightarrow } \varphi_i,
$$
where the vector $\varphi$ belongs to $\mathcal F$ and  is characterized by
$$
\log (\rho) = R(\varphi) - I(\varphi).
$$
The occupancy frequency $\varphi$ coincides with the stationary distribution $u$ of  $X$ only if $m_{\jmath(P)}=m_{\jmath(P')}$ for all patches $P,P'$.
\end{thm}

The results obtained for the periodic and random environments can be derived similarly to the infinite transitive graphs.

\appendix

\section{Proof of Theorem \ref{rateped}}
\label{proof}
Recall that $(X_n)$  is an irreducible, aperiodic Markov chain with transition matrix $D$. As the state space $\{1,\ldots, K\}$ is finite, $(X_n)$ is strongly irreducible and there exists $n_0\in \NN$ such that $D^{n_0}>>0$ (all positive entries).

Now the real number $I(f)$ is the cost for the habitat occupation frequencies associated with the random walk $X$ to equal $f$. It gives the geometric decrease
of the probability that the portion of time  spent in habitat $i$ until generation $n$ is close to $f_i$ :
\begin{eqnarray}
\label{frqc}
I(f_1,\ldots ,f_K) = \lim_{\epsilon \rightarrow 0} \lim_{n\rightarrow \infty} -\frac{1}{n} \log \PP(   f_k-\epsilon \leq F_k(X_n) \leq f_k+\epsilon).
\end{eqnarray}
 This result holds for any irreducible aperiodic random walk $X$ when $I$ is finite. Indeed, thanks to Sanov's theorem (see e.g. \cite[Theorem 3.1.6 page 62]{dembo}, $I$ is convex and continuous.  The function $I$ is called the rate function associated to the path of the random walk $X$.\\

Finally, the fact that the offspring distribution $N_i$ for an individual living  in patch $i$ satisfies $\EE(N_i\log^+ N_i)<\infty$ ensures that for all $i,j=1,\ldots, K$,
$$\EE( Z_{1}^{j}\log^+  Z_1^{(j)} \ \vert \ Z_0^{(i)}=1, \ Z_0^{(k)}=0 \text{ for } k\ne i )<\infty .
$$
Then $Z_n/\rho^n$ converges to a non degenerate variable, see e.g. \cite[Chapter 5, section 6, Theorem 1]{AN}.

\paragraph{Expression of  the growth rate $\rho$ and  habitat occupation frequencies.}
 We start with one individual in patch $1$. 
Recalling that for every $i=1,\ldots ,K$,  $\EE( Z_{n}^{(i)}  )=\sum_{j=1}^K\EE( Z_{n-1}^{(j)})m_{j}d_{ji}$, we get by induction
$$\EE( Z_{n}^{(i)} )=\sum_{\substack{j_0=1, \ j_n=i, \\ \quad  1\leq j_1,\ldots ,j_{n-1} \leq K} } 
\Pi_{k=0}^{n-1}m_{j_k}d_{j_k j_{k+1}}=\EE_1\big(1_{X_n=i}\Pi_{i=0}^{n-1}m_{X_i}\big).$$
This yields
$$\EE(\vert Z_n \vert )= \EE\left(\Pi_{i=0}^{n-1} m_{X_i}\right)$$
Denoting by $S_n(i)$ the number of  visits  of the random disperser $X$ in patch $i=1,\ldots V,$ :
$$S_n(i)=\#\{ k\leq n-1 : \jmath (X_k)=i\}=nF_i(X_n),$$
we deduce
\begin{eqnarray*}
\EE(\vert Z_n \vert )
&=&\EE\left(\Pi_{i=1}^K m_{i}^{S_n(i)}\right) \\
&=&\int_{\mathcal{F}} \exp(n\sum_{i=1}^Kf_i\log(m_i))\PP(F_1(X_n) \in df_1,\ldots ,F_K(X_n)\in df_K). 
\end{eqnarray*}
Using $(\ref{frqc})$ and Laplace method, we get  
$$\log(\EE(\vert Z_n \vert)^{1/n}) \stackrel{n\rightarrow\infty}{\longrightarrow}  \max\{ \sum_{i=1}^K f_i\log(m_i) -I(f_1,\ldots ,f_K) : f\in \mathcal F\}.$$
This proves the first part of the result.\\

The  maximum of $h:=R-I$ is reached for a unique frequency $\varphi$, which means that there is a unique $\varphi$ such that $\log(\rho)=R(\varphi)-I(\varphi)$. This is the object of Proposition \ref{uniq}. Moreover  the partial derivatives of $h$ at $\varphi$ are zero. As $f_K=1-f_1-\ldots f_{n-1}$, for every $1\leq i\leq K-1$,
$$
\log(m_i)-\log(m_K)  -\frac{\partial}{\partial f_i} I_{\vert f=\varphi}=0.
$$
If $\varphi=p$, then the partial derivatives of $I$ at $\varphi$ are zero. This can be directly computed or deduced from   $I\geq 0$ and $I(p)=0$. This ensures that for every $1\leq i\leq K-1$, 
$\log(m_i)-\log(m_K)=0$, i.e. $m_1=m_2=\ldots =m_{K}$. This proves the third part of the theorem.

Finally, let us prove that the habitat occupation frequencies  of a typical individual are given by the vector $\varphi$. Let $\epsilon >0$ and $i \in\{1,\ldots,K\}$.
The individuals alive in generation $n$ are labeled by $\mathbf{u}_k$, $k=1,\ldots ,Z_n$.

Following the first part of the proof,
\begin{eqnarray*}
\EE(\sum_{k=1}^{Z_n} 1_{\vert F_i(\mathbf{u}_k)-\varphi_i \vert\geq \epsilon})&=& \EE\left(1_{\vert F_i(X_n)-\varphi_i \vert\geq \epsilon}\Pi_{i=0}^n m_{X_i} \right)\\
&=&\EE\left(1_{\vert F_i(X_n)-\varphi_i \vert\geq \epsilon}\Pi_{i=1}^K m_{i}^{S_n(i)}\right) \\
&=&\int_{\mathcal{F}} 1_{\vert f_i-\varphi_i \vert\geq \epsilon}\exp(n\sum_{i=1}^Kf_i\log(m_i))\PP(F_1(X_n) \in df_1,\ldots ,F_K(X_n)\in df_K).
\end{eqnarray*}
Using again  $(\ref{frqc})$ and the Laplace method, we get  
$$\frac{1}{n}\log \EE\big(\sum_{k=1}^{Z_n}
1_{\vert F_i(\mathbf{u}_k)-\varphi_i \vert\geq \epsilon}\big) \stackrel{n\rightarrow\infty}{\longrightarrow} C_{i,\epsilon},$$
with $C_{i,\epsilon}= \max\{ \sum_{i=1}^K f_i\log(m_i) -I(f_1,\ldots ,f_K) : f\in \mathcal F, \ \vert f_i-\varphi_i \vert\geq \epsilon\}.$
The  uniqueness of the argmax $\varphi$ ensures that 
$C_{i,\epsilon}< C_{i,0}$. Moreover the growth rate $C_{i,0}$ is equal to $\log(\rho)$ and
$$\frac{1}{n}\log \EE\big(\sum_{k=1}^{Z_n}
1_{\vert F_i(\mathbf{u}_k)-\varphi_i \vert\geq \epsilon}\big)-\frac{1}{n}\log \rho^n\stackrel{n\rightarrow\infty}{\longrightarrow} C_{i,\epsilon}-C_{i,0}<0.$$
Then 
$$\EE\big(  \frac{1}{\rho^n}.\sum_{k=1}^{Z_n} 1_{\vert F_i(\mathbf{u}_k)-\varphi_i \vert\geq \epsilon}\big)\stackrel{n\rightarrow\infty}{\longrightarrow}0.$$
In other words $\sum_{k=1}^{Z_n} 1_{\vert F_i(\mathbf{u}_k)-\varphi_i \vert\geq \epsilon}/\rho^n$ goes to $0$ in probability.
Adding that $Z_n \sim W \rho^n$ a.s. as $n\rightarrow \infty$ and $\{W>0\}=\{\forall n \in \NN, \  Z_n>0\}$ a.s. ensures that 
 $$ 1_{Z_n>0} \frac{1}{Z_n}.\sum_{k=1}^{Z_n} 1_{\vert F_i(\mathbf{u}_k)-\varphi_i \vert\geq \epsilon}\stackrel{n\rightarrow\infty}{\longrightarrow}0$$
in probability.  By dominated convergence,
$$ \EE\left(1_{Z_n>0} \frac{1}{Z_n}.\sum_{k=1}^{Z_n} 1_{\vert F_i(\mathbf{u}_k)-\varphi_i \vert\geq \epsilon} \right)\stackrel{n\rightarrow\infty}{\longrightarrow}0.$$
Then,  conditionally on $Z_n>0$, denoting  by $\mathbf{U}_n$ an individual chosen uniformly in generation $n$,
$$\PP\left(\vert F_i(\mathbf{U}_n)-\varphi_i \vert\geq \epsilon, \quad Z_n>0\right)\stackrel{n\rightarrow\infty}{\longrightarrow}0.$$
This proves that $ F_i(\mathbf{U}_n)\stackrel{n\rightarrow \infty}{\longrightarrow } \varphi_i$ in probability
and completes the proof.

\paragraph{Study of $I$ and uniqueness of argmax $R-I$.}
The supremum $I$ defined by
$$
I(f)=I(f_1,\ldots ,f_K):=\sup\{\sum_{j=1}^K f_j\log(u_j/(uD)_j) : u\in \RR^K, u>>0\},$$
is reached for a unique unit positive vector. This means that there exists a unique $u(f)=(u_1,\ldots, u_K)$ such that
 $$I(f)=\sum_{j=1}^K f_j\log(u_j(f)/(u(f)D)_j), \quad u_1(f)+\cdots+u_K(f)=1, \quad  u_1(f)>0, \ldots, u_K(f)>0.$$ Indeed
this vector $u(f)$ realizes  a maximum  for $u\in \RR^K, u>>0$ and thus satisfies for $j=1,\ldots ,K$,
\begin{equation} \label{idt}
 \frac{f_j}{u_j}-\sum_{i=1}^{K}d_{ji} \frac{f_i}{(uD)_i}=0.
\end{equation}
This equation  characterizes $u$, see  Exercise IV.9 page 46 in \cite{DenH}, which ensures that $u(f)$ is uniquely defined. Note also that if $f_i$ is the stationary distribution, $u_j=f_j$ satisfies this equation since $(fD)_i=f_i$, so that $I=0$. \\

\begin{prop}\label{uniq}
There exists a unique $\varphi \in \mathcal{F}$ such that $\log(\rho)=R(\varphi)-I(\varphi)$.
\end{prop}
\paragraph{Proof.} 
We observe that $f\mapsto u(f)$ can be extended from $\mathcal{F}$ to $[0,\infty)^K\setminus\{0\}$ and satisfy $(\ref{idt})$
on $[0,\infty)^K\setminus\{0\}$ by setting 
$$u(f)=u(f/ \parallel f \parallel ), \quad{where}   \ \ \parallel f \parallel = \sum_{i=1}^K f_i.$$ 
Then  $R-I$ can be extended   to  $[0,\infty)^K\setminus\{0\}$ with
$$(R-I)(f)= \sum_{j=1}^K \frac{f_j}{\sum_{k=1}^K f_k}\log(m_j(u(f)D)_j/u_j(f)),$$
and $(R-I)(\lambda f)=(R-I)(f)$ for every $\lambda \in (0,\infty)$.

Consider a vector $f$ which realizes the maximum of $R-I$ and does not belong to the boundary of $[0,\infty]^K$. Then the partial derivatives are zero 
 and for every $i=1, \ldots, K$,
$$\frac{1}{\sum_{k=1}^K f_k}.
\bigg[\log(m_i(uD)_i/u_i) +\sum_{j=1}^K f_j
\left(
\frac{\frac{\partial}{\partial f_i} (uD)_{j}}{(uD)_j}-\frac{\frac{\partial}{\partial f_i}u_j}{u_j}
\right)
\bigg]=\frac{\sum_{j=1}^K f_j\log(m_j(uD)_j/u_j)}{[\sum_{k=1}^K f_k]^2}.$$
Using (\ref{idt}) we get
$$
\sum_{j=1}^K f_j
\frac{\frac{\partial}{\partial f_i} (uD)_{j}}{(uD)_j}
=\sum_{j=1}^{K} \frac{f_j}{(uD)_j}\sum_{k=1}^K d_{kj} \frac{\partial u_k}{\partial f_i} 
= \sum_{k=1}^K   \frac{\partial u_k}{\partial f_i}  \sum_{j=1}^{K} \frac{f_jd_{kj}}{(uD)_j} 
=\sum_{k=1}^{K} \frac{\partial u_k}{\partial f_i}\frac{f_k}{u_k},
$$
so that
$$\log(m_i(uD)_i/u_i) =\frac{\sum_{j=1}^K f_j\log(m_j(uD)_j/u_j)}{\sum_{k=1}^K f_k}.$$
Observe that the right hand side does not depend on $i$, so that for every $i=1, \ldots, K$,
$$ (uD')_i=\alpha u_i,$$
where $D'_{ji}=d_{ji}m_i$ and $\alpha$ is a positive constant. Then $u$ is left eigenvector of $D'$ with positive entries.
Moreover $D'$ is strongly irreducible since $D$ is strongly irreducible \footnote{$D$ is strongly irreducible if there exists $n_0\geq 1$ such that all the coefficients of $D^{n_0}$
are positive. This property is a consequence of irreducibility and aperiodicity.}
 and $m_i>0$ for every $i$ by assumption.
Now Perron--Frobenius theory  ensures that there is a \emph{unique} left positive eigenvector  $u$
of $D'$ such that $\sum_{i=1}^K u_i=1$ (see e.g. \cite{HornJ} \footnote{which states that two positive eigenvectors of a primitive matrix are colinear. It comes actually from the classical decomposition of $A^n$
using the maximum eigenvalue and the  left and right eigenvector associated with.}.  Moreover, following the literature on large deviations (see e.g. \cite{DenH}), (\ref{idt}) reads
$$
f_j=\sum_{i=1}^{K}D''_{ji} f_i
$$
with $D''_{ji}=u_jd_{ji}/(uD)_i$.  Here again $D''$ is strongly irreducible since $D$ is strongly irreducible and both $u$ and $uD$ are positive vectors. Using again Perron--Frobenius theory guarantees the  uniqueness of the solution
$f$ such that $\sum_{i=1}^K f_i=1$. This ensures the uniqueness of the argmax of $R-I$ in the interior of $[0,\infty)^K$. We complete the proof by adding that there is at least one argmax in the interior of $[0,\infty)^K$ since we recall that the frequency occupation is  the product of the right and left eigenvectors of $A$ associate to $\rho$, which are both positive.
If $\varphi_1$ and $\varphi_2$ realize the max of $R-I$, the concavity of this function ensures that the segment $[\varphi_1, \varphi_2]$, which is contradiction with the uniqueness  in the interior of $[0,\infty)^K$ and complete the proof.
 \hfill$\Box$

\section{Annex : Classification Theorem for MBPRE}
We consider here a multitype branching process in random environment $Z_n=(Z_n^{(i)} : i=1, \ldots,K)$ whose mean offspring matrix is denoted by $A=A(w)$, where $w=(w_0,w_1,\ldots)$ is the environment.\\
We introduce the extinction probability vector in environment $w$  starting from one individual in habitat $i$:
$$q_i(w)=\lim_{n\rightarrow \infty} \PP( \vert Z_n \vert =0 \ \vert \ w, \  Z_0^{(i)}=1, \  Z_0^{j}=0 \ \text{if } j\ne i).$$
\begin{prop}[\cite{Tanny}, Theorems 9.6 and 9.10] Assuming that
$$\PP( q_i(w)<1 : i =1,\ldots, K )= 1 \ \ \text{or} \ \ \PP( q_i(w) = 1 :  i =1,\ldots, K)= 1, \qquad (*)$$
we have   
\begin{itemize}
\item If $\gamma<0$, then  the probability of extinction is equal to $1$ for almost every $w$.
\item If $\gamma=1$, then
\begin{itemize}
\item either  for every $m \geq 1$, 
$w$-a.s., there exists  $1\leq i\leq K$ such that  

$\PP(\vert Z_m \vert > 1  \ \vert \ \ w, \ Z_0^{(i)}=1, \  Z_0^{(j)}=0 \ \text{if } j\ne i)=0$.
\item or  $q_i(w)=1$ $w$-a.s..
\end{itemize} 
\end{itemize}
Assuming that there exist integers $N,L>0$   such that $\PP( \forall 1\leq i,j\leq K, \ (A_N\cdots A_0)_{ij}\not=0   ) = 1$ and
$\vert \EE (\log (1-\PP(Z_K = 0 \ \vert \ Z_0^{(L)}=1) ) \vert <\infty$,  then $(*)$ is satisfied and
\begin{itemize}
\item If $\gamma>0$,   then $w$-a.s $q_i(w)<1$ for every $i=1,\cdots,K$, and 
$$\PP\big(\lim_{n\rightarrow \infty} n^{-1}\log(\vert Z_n\vert)=\gamma \ \vert  \ w,  \ Z_0^{(i)}=1, \  Z_0^{j}=0 \ \text{if } j\ne i \big)=1-q_i(w).$$
\end{itemize} 
\end{prop}
Moreover thanks to \cite[Theorem 9.11]{Tanny}, if all the coefficients of the matrix $A$ 
are positive and bounded, i.e.,
$$\exists 0<c,c'<\infty, \ c\leq \inf_{1\leq i,j\leq K} A_{i,j}\leq \sup_{1\leq i,j\leq K} A_{ij}\leq <c' \quad \text{a.s.},$$
then
$$Z_n= O(\| A_{n-1}\cdots A_0 \|)  \quad  \text{a.s.}$$

\textbf{Acknowledgement.} This work was funded by project MANEGE `Mod\`eles
Al\'eatoires en \'Ecologie, G\'en\'etique et \'Evolution'
09-BLAN-0215 of ANR (French national research agency).

\end{document}